\documentclass[twocolumn]{cifa}
\usepackage[latin1]{inputenc}
\usepackage[T1]{fontenc}
\usepackage[francais]{babel}
\usepackage{graphicx}
\usepackage{psfrag}
\usepackage{color}
\usepackage{epsfig}
\usepackage[hang]{caption2}
\usepackage[normalsize,it,hang]{subfigure}
\usepackage{amsmath}
\usepackage{array}
\usepackage{theorem}
\usepackage{amsmath,amssymb}
\usepackage{amsfonts}
\usepackage{amssymb}
\usepackage{subfigure}
\usepackage{amssymb,amsbsy,amsmath,amsfonts,amssymb,amscd}

\newtheorem{remarque}{\it Remarque\/}

\begin{document}

\title{Vers une commande sans mod\`{e}le \\ pour am\'{e}nagements hydro\'{e}lectriques en
cascade}

\author{\vskip 1em
C\'{e}dric \textsc{JOIN}\textsuperscript{1},
G\'{e}rard \textsc{ROBERT}\textsuperscript{2},
Michel \textsc{FLIESS}\textsuperscript{3}\\
\vskip1em
\textsuperscript{1}INRIA-ALIEN \& CRAN (UMR CNRS 7039)\\
Nancy-Universit\'{e}, BP 239, 54506 Vand\oe{}uvre-l\`es-Nancy, France \\
\textsuperscript{2}EDF, Centre d'Ing\'{e}ni{e}rie Hydraulique \\ Savoie
Technolac, 73373 Le Bourget du Lac, France\\
\textsuperscript{3} INRIA-ALIEN \&  LIX
(UMR CNRS 7161) \\ \'Ecole polytechnique, 91128 Palaiseau, France  \\
\vskip 1em
{ \texttt{Cedric.Join@cran.uhp-nancy.fr, gerard.robert@edf.fr,
\\ Michel.Fliess@polytechnique.edu\\}}
} \maketitle

\begin{abstract}
On aborde la r\'{e}gulation du niveau d'eau dans un am\'{e}nagement
hydraulique, soumis \`{a} de fortes contraintes,  par la {\em commande
sans mod\`{e}le}. Les nombreuses simulations num\'{e}riques fournissent
d'excellents r\'{e}sultats, obtenus gr\^{a}ce \`{a} des algorithmes robustes et
simples.
\\
{\it Abstract}---A new concept called {\it Model-Free Control} is
applied to hydroelectric run-of-the river power plants, with severe
constraints and operating conditions. Numerous computer simulations
display excellent results, which are obtained thanks to simple and
robust algorithms.
\end{abstract}


\section{Introduction}

La r\'{e}gulation du niveau d'eau dans un am\'{e}nagement hydraulique
(barrages, retenues, rivi\`{e}res, canaux, \dots), importante pour la
production d'\'{e}lectricit\'{e}, la navigation, l'irrigation, et bien
d'autres usages, a re\c{c}u une attention consid\'{e}rable, selon des points
de vue les plus vari\'{e}s, dans la litt\'{e}rature automaticienne r\'{e}cente:
voir, par exemple,
\cite{bastin,bes,cantoni,coron,dang,dos,dul,dumur,georges0,halleux,ham,lem,rab,rab0,thomassin},
et les r\'{e}cents tours d'horizon \cite{georges,litrico,mar} et
\cite{zhuan}. La difficult\'{e} tient \`{a} la nature m\^{e}me du proc\'{e}d\'{e}:
\begin{itemize}
\item c'est un syst\`{e}me hydraulique d'\'{e}coulement \`{a} surface libre, encha\^{\i}n\'{e} (voir
Fig. \ref{bief}) ou non, avec g\'{e}om\'{e}trie quelconque du volume d'eau;
\begin{figure*}
\center{\includegraphics[width=1.575\columnwidth]{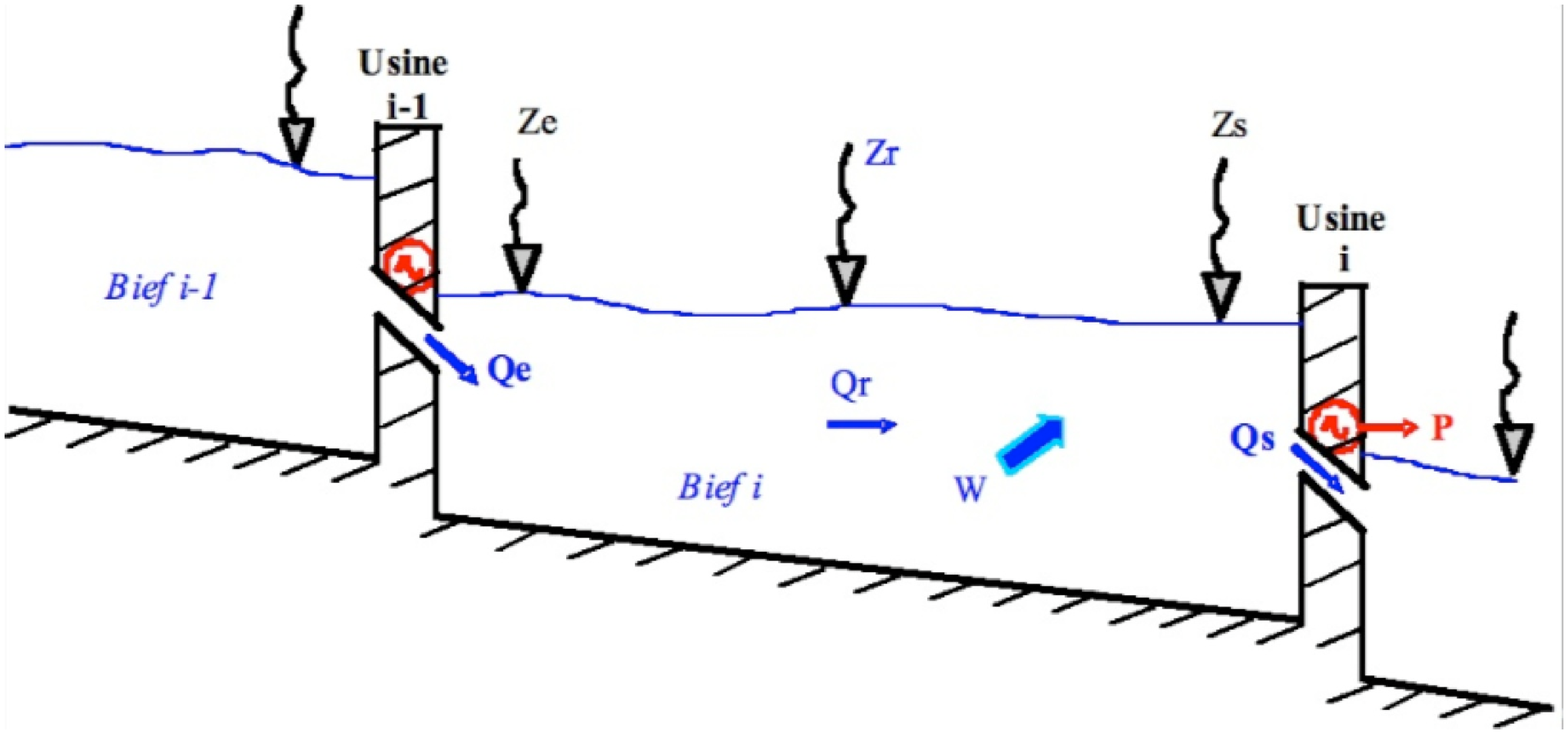}}
 \caption{Sch\'{e}ma d'un bief\label{bief}}
\end{figure*}
\item sa taille est grande (plusieurs kilom\`{e}tres de long);
\item il est non lin\'{e}aire, multivariable et distribu\'{e} dans
l'espace (\'{e}quation aux d\'{e}riv\'{e}es partielles de Saint-Venant);
\item il existes des perturbations al\'{e}atoires dont l'amplitude peut \^{e}tre
importante, mais inconnue, dues \`{a} des apports ou pr\'{e}l\`{e}vements d'eau.
\end{itemize}
L'embarras cro\^{\i}t encore si
\begin{itemize}
\item le niveau d'eau \`{a} r\'{e}guler est \'{e}loign\'{e} de l'organe r\'{e}glant
(une vanne, par exemple), compte-tenu des retards variables mais
aussi de ph\'{e}nom\`{e}nes hydrauliques complexes de basculement des plans
d'eau, en fonction du d\'{e}bit;

\item le cahier des charges exige une
robustesse des performances sur toute la plage de fonctionnement et
vis-\`{a}-vis des variations des caract\'{e}ristiques de l'am\'{e}nagement
(temps de retard variables, superficie du r\'{e}servoir variable en
fonction du niveau et de l'envasement, usure des organes r\'{e}glants,
\dots).
\end{itemize}
Cette communication, qui compl\`{e}te une pr\'{e}c\'{e}dente \cite{portugal},
traitant d'un cas simplifi\'{e}, propose pour des am\'{e}nagements
hydro\'{e}lectriques en cascade, s\'{e}par\'{e}s par des biefs (voir Fig.
\ref{bief}), la nouvelle \og commande sans mod\`{e}le \fg ~
(\cite{esta,malo}), facile \`{a} mettre en {\oe}uvre, et \'{e}minemment robuste,
dont nous rappelons l'essentiel au {\S} \ref{rappel}. Apr\`{e}s une
description du dispositif et de sa r\'{e}gulation au {\S}
\ref{description}, le {\S} \ref{simu} fournit d'excellentes simulations
num\'{e}riques, respectant plusieurs sc\'{e}narios. Le {\S} \ref{conclusion}
conclut au grand int\'{e}r\^{e}t pratique de notre approche, tout en posant
la question du r\^{o}le de la dimension infinie.

\vspace{0.1cm} {\small \noindent{\bf Remerciements}. Travail sous
l'\'{e}gide de deux contrats, intitul\'{e}s \og Commande sans mod\`{e}le pour
am\'{e}nagements hydro\'{e}lectriques en cascade \fg, entre, d'une part, EDF
et, d'autre part, l'\'Ecole polytechnique, l'INRIA et le CNRS. Le
syst\`{e}me de r\'{e}gulation, qui y est d\'{e}crit, a fait l'objet d'une
demande de brevet (n° FR0858532), d\'{e}pos\'{e}e par EDF et l'\'Ecole
polytechnique le 12 d\'{e}cembre 2008.}

\section{Rappels}\label{rappel}
\subsection{Commande sans mod\`{e}le et i-PID}
La {\em commande sans mod\`{e}le} repose sur une mod\'{e}lisation {\em
locale}, sans cesse r\'{e}actualis\'{e}e, \`{a} partir de la seule connaissance
du comportement entr\'{e}e-sortie. \`A l'\'{e}quation diff\'{e}rentielle
inconnue, lin\'{e}aire ou non,
\begin{equation}\label{E}
E (y, \dot{y}, \dots, y^{(a)}, u, \dot{u}, \dots, u^{(b)}) = 0
\end{equation}
d\'{e}crivant approximativement le comportement entr\'{e}e-sortie, on
substitue le mod\`{e}le \og ph\'{e}nom\'{e}nologique \fg, valable sur un court
laps de temps,
\begin{equation}\label{F}
y^{(\nu)} = F + \beta u
\end{equation}
o\`{u}
\begin{itemize}
\item l'ordre de d\'{e}rivation $\nu$, en g\'{e}n\'{e}ral $1$ ou $2$, donc diff\'{e}rent
de l'ordre de d\'{e}rivation $a$ de $y$ en (\ref{E}), est fix\'{e} par
l'op\'{e}rateur;
\item le param\`{e}tre constant $\beta$ est fix\'{e} par
l'op\'{e}rateur de sorte qu'en (\ref{E}) $\beta u$ et $y^{(\nu)}$ aient
m\^{e}me ordre de grandeur.
\end{itemize}
La valeur de $F$ \`{a} chaque instant se d\'{e}duit de celles de $u$ et de
$y^{(\nu)}$, obtenus par d\'{e}rivateurs num\'{e}riques. On obtient le
comportement d\'{e}sir\'{e}, si, par exemple, $\nu = 2$ en \eqref{F}, gr\^{a}ce
au {\em correcteur} {\em PID intelligent},  ou, en abr\'{e}g\'{e}, {\em
i-PID},
\begin{equation}\label{universal}
u = - \frac{F}{\beta} + \frac{\ddot{y}^\ast}{\beta} + K_P e + K_I
\int e + K_D \dot{e}
\end{equation}
o\`{u}
\begin{itemize}
\item $y^\ast$ est la trajectoire de r\'{e}f\'{e}rence de la sortie, obtenue
selon les pr\'{e}ceptes de la commande par platitude;
\item $e = y - y^\ast$ est l'erreur de poursuite;
\item $K_P$, $K_I$, $K_D$ sont les gains de
r\'{e}glage.
\end{itemize}
\begin{remarque}\label{trivial}
Il est ais\'{e} de d\'{e}terminer ces gains car, avec \eqref{universal},
\eqref{F} se ram\`{e}ne, contrairement aux PID classiques, \`{a} un
int\'{e}grateur pur du second ordre.
\end{remarque}
\begin{remarque}
Voir
\begin{itemize}
\item \cite{esta,malo} pour plus de d\'{e}tails,
\item \cite{med} pour les liens avec les PID usuels,
\item \cite{mines,choi,brest1,brest2,psa,vil} pour plusieurs autres
illustrations concr\`{e}tes.
\end{itemize}
\end{remarque}

\subsection{D\'{e}rivateurs num\'{e}riques}

\subsubsection{Un calcul
simple}\label{calcul} Avec les notations classiques du calcul
op\'{e}rationnel ({\it cf.} \cite{yosida}), il correspond \`{a}
$p(t)=a_{0}+a_{1}t$, $a_0, a_1 \in \mathbb{R}$, pour $t \geq 0$:
\begin{equation}\label{eq5}
P(s)=\frac{a_{0}}{s}+\frac{a_{1}}{s^{2}}
\end{equation}
On cherche \`{a} \'{e}liminer $a_{0}$ car on veut estimer $a_{1}$. Pour cela
on multiplie \eqref{eq5} par $s$:
    \begin{equation*}\label{eq6}
      sP(s)=a_{0}+\frac{a_{1}}{s}
    \end{equation*}
puis on d\'{e}rive l'expression obtenue par rapport \`{a} $s$ pour supprimer
$a_0$:
    \begin{equation*}\label{eq7}
      P(s)+s\frac{\text{d}P(s)}{\text{d}s}=-\frac{a_{1}}{s^{2}}.
    \end{equation*}
Avant de revenir au domaine temporel, une multiplication par
$s^{-N}$, avec $N>1$, $N=2$ par exemple, est n\'{e}cessaire pour obtenir
uniquement des int\'{e}grales:
    \begin{equation*}\label{eq8}
      s^{-2}P(s)+s^{-1}\frac{\text{d}P(s)}{\text{d}s}=-s^{-4}a_{1}.
    \end{equation*}
On revient au domaine temporel en rappelant ({\it cf.}
\cite{yosida}) que $\frac{d}{ds}$ correspond \`{a} la multiplication par
$-t$:
\begin{equation}
\begin{split}
        a_{1}&=\frac{6\Bigg(\displaystyle\int_{t_{0}}^{t}\tau x(\tau)\text{d}\tau-\int_{t_{0}}^{t}\int_{t_{0}}^{\tau}x(\kappa)\text{d}\kappa\text{d}\tau\Bigg)}{\displaystyle t^{3}}\\&=\frac{6\Bigg(\displaystyle\int_{t_{0}}^{t}\tau x(\tau)\text{d}\tau-\int_{t_{0}}^{t}(t-\tau)x(\tau)\text{d}\tau\Bigg)}{\displaystyle t^{3}}\\&=\frac{6\displaystyle\int_{t_{0}}^{t}\big(\tau x(\tau)-(t-\tau)x(\tau)\big)\text{d}\tau}{\displaystyle t^{3}}
\end{split}\label{eq15}
\end{equation}

\subsubsection{Filtres d\'{e}rivateurs} Les propri\'{e}t\'{e}s du calcul
op\'{e}rationnel permettent de g\'{e}n\'{e}raliser \eqref{eq15} en \'{e}crivant tout
estimateur alg\'{e}brique d'une d\'{e}riv\'{e}e d'ordre quelconque d'un signal
$x(t)$, analytique autour de $0$, sous la forme d'un filtre
$\int_{t_{1}}^{t_{n}} \varpi(t) x(t)\text{d}t$, o\`{u} $\varpi(t)$ est
un polyn\^{o}me temporel caract\'{e}risant ledit estimateur. Les calculs du
{\S} \ref{calcul} permettent de comprendre la vari\'{e}t\'{e} possible
d'estimateurs que l'on peut obtenir par cette m\'{e}thode, en modifiant,
par exemple, l'ordre de troncature du d\'{e}veloppement de Taylor.

\begin{remarque}
Pour compl\'{e}ter ce rapide r\'{e}sum\'{e}, inspir\'{e} de \cite{gretsi}, renvoyons
\`{a} \cite{mboup} pour des d\'{e}veloppements th\'{e}oriques importants, y
compris pour l'implantation algorithmique. Voir \cite{easy} pour
plus de r\'{e}f\'{e}rences et, surtout, pour les progr\`{e}s substantiels que
permet cette d\'{e}rivation en automatique non lin\'{e}aire, avec mod\`{e}le
connu, aux incertitudes param\'{e}triques pr\`{e}s.
\end{remarque}

\section{Description du dispositif}\label{description}
\subsection{G\'{e}n\'{e}ralit\'{e}s}
\begin{figure}
\center{\includegraphics[width=1.08\columnwidth]{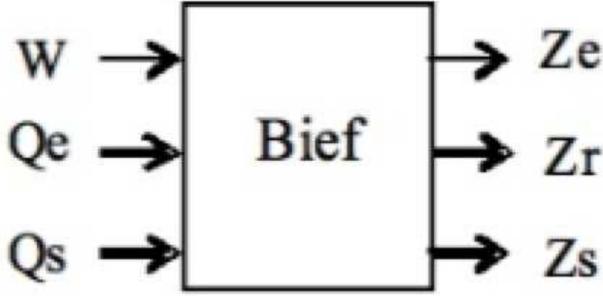}}
 \caption{Sch\'{e}ma bloc du mod\`{e}le \label{model}}
\end{figure}

Comme l'indiquent les sch\'{e}mas \ref{bief} et \ref{model}, on cherche
\`{a} ma\^{\i}triser le niveau d'eau $z_r$ du bief $i$, entre les usines
$i-1$ et $i$, avec
\begin{itemize}
\item commande $u$, ou d\'{e}bit de sortie $Q_s$, \'{e}chantillonn\'{e}e-bloqu\'{e}e
et satur\'{e}e en position et vitesse;
\item deux entr\'{e}es exog\`{e}nes $Q_e$ et $W$, perturbatrices;
\item insensibilit\'{e} du capteur de niveau qui a une r\'{e}solution de 1
cm.
\end{itemize}
La Fig. \ref{cde} fournit le sch\'{e}ma de commande:
\begin{figure}[h!h]
\center{\includegraphics[width=1.03\columnwidth]{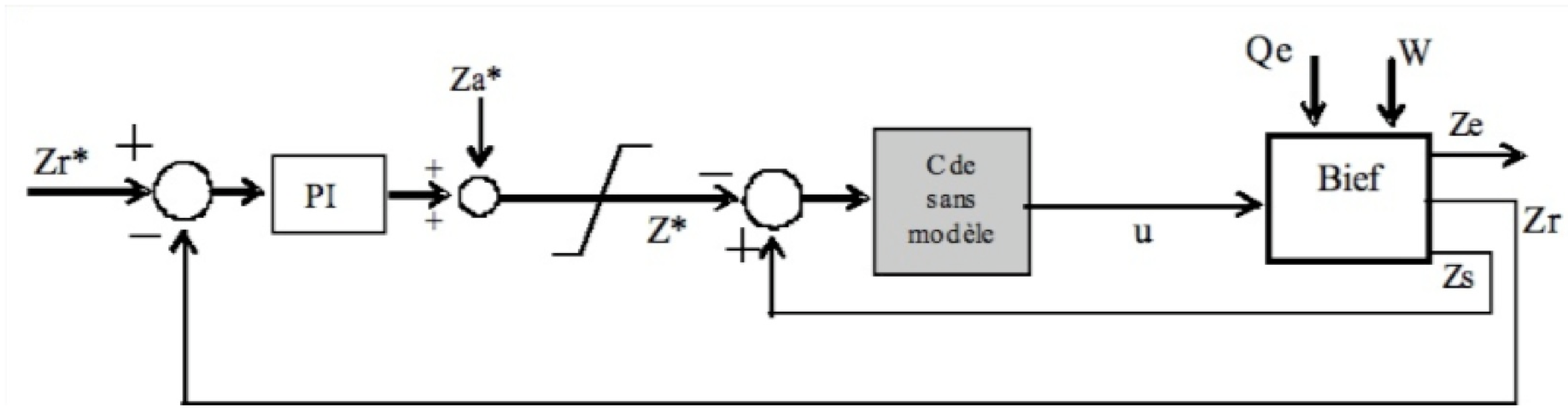}}
 \caption{Sch\'{e}ma bloc de la commande d\'{e}sir\'{e}e \label{cde}}
\end{figure}

\subsection{Traitement du retard}
Le d\'{e}bit $Q_e$ modifie le niveau $z_r$, loin de l'usine $i$, avec un
retard mal connu et, qui plus est, variable, donc difficile, sinon
impossible, \`{a} ma\^{\i}triser avec les outils th\'{e}oriques actuels
d'identification et de commande. On y pallie en reconstruisant un
niveau $z_a$ \`{a} partir de mesures proches de l'usine $i-1$, par deux lois
possibles:
\begin{enumerate}
\item l'une repose sur la connaissance de niveaux
\begin{equation*}z_a=f(Z_e)\label{zn}\end{equation*}
\item l'autre sur celle du d\'{e}bit entrant
\begin{equation*}z_a=g(Q_e)\label{zq}\end{equation*}
\end{enumerate}
fournies par EDF et qui r\'{e}sultent de lois empiriques.

La trajectoire de r\'{e}f\'{e}rence, interm\'{e}diaire, $z^\star_a$, d\'{e}termin\'{e}e
selon les principes de la platitude, est compl\'{e}t\'{e}e par un correcteur
pilot\'{e} par l'erreur de poursuite sur $z_r$. Nous r\'{e}alisons, ainsi,
une \textit{replanification} en ligne automatique de la trajectoire
de r\'{e}f\'{e}rence. Cette trajectoire corrig\'{e}e $z^\star$ devient la
consigne d'une boucle interne, \`{a} commande sans mod\`{e}le.

\subsection{Commande sans mod\`{e}le}
Nous mettons ici en {\oe}uvre \eqref{F}, avec $\nu = 1$, sous la forme
\begin{equation}
\alpha \dot y = F-u
\label{sm}
\end{equation}
qui a l'avantage de donner un sens clair de d\'{e}bit entrant \`{a} $F$,
information importante pour l'exploitation de la centrale
hydraulique, en faisant appara\^{\i}tre la diff\'{e}rence $F-u$, o\`{u} $u$ est
le d\'{e}bit sortant. La boucle est ferm\'{e}e, d'apr\`{e}s \eqref{universal},
par
$$
u = F -\alpha\dot{y}^\star + \alpha \mbox{\rm PI}(e)
$$
o\`{u} $\mbox{\rm PI}(e) = K_P e + K_I\int e$ est un correcteur PI,
pilot\'{e} par $e = y-y^\star$.

\section{Quelques r\'{e}sultats de simulation}\label{simu}

\subsection{Sc\'{e}narios}
Trois sc\'{e}narios, propos\'{e}s par EDF sur l'am\'{e}nagement hydro\'{e}lectrique
de Fessenheim, permettent d'\'{e}valuer les performances de notre
commande. Cet am\'{e}nagement au fil de l'eau (capable de turbiner
$175,5$ MW en pointe) fait partie des 10 centrales EDF (avec sa
filiale EnBW) encha\^{\i}n\'{e}es en cascade sur le Rhin et r\'{e}gul\'{e}es en
niveau. Il est caract\'{e}ris\'{e} par un bief de 15 kilom\`{e}tres de long et
d'une hauteur de chute d'environ 15 m\`{e}tres. Le principal objectif \`{a}
satisfaire est le respect de la contrainte
\begin{equation}\label{ineq}
z_r^\star-10 ~\mbox{\rm cm} < z_r < z_r^\star+10 ~\mbox{\rm cm}
\end{equation}
en un point distant de 7,5 kilom\`{e}tres de l'actionneur.
\begin{enumerate}
\item Le premier sc\'{e}nario repr\'{e}sente la fin d'une crue avec de fortes
variations de d\'{e}bit.
\item Le second, plus doux, correspond \`{a} une situation normale plus fr\'{e}quente;
\item Quant au troisi\`{e}me, il est acad\'{e}mique afin de mettre en
lumi\`{e}re les comportements lors de saturations.
\end{enumerate}
Les sc\'{e}narios 1 et 2, qui sont tout \`{a} fait r\'{e}alistes, ont une dur\'{e}e
de $4$ jours. Les figures peuvent donner l'illusion, donc, d'une
commande tr\`{e}s dynamique, ce qui n'est pas vraiment le cas.

Les perturbations sont constitu\'{e}es, d'une part, par un biais \'{e}gal \`{a}
$0.03 Q_e+10$ pour prendre en compte les erreurs de d\'{e}bit (consigne
- mesure) et, d'autre part, par des sass\'{e}es de 100 $\mbox{\rm
m}^3/\mbox{\rm s}$, reproduisant le comportement d'une \'{e}cluse. Ces
derni\`{e}res sont tr\`{e}s violentes puisqu'elles ne durent que $15$
minutes, c'est-\`{a}-dire $7$ \'{e}chantillons.

\subsection{R\'{e}glage de la boucle externe}
Nous avons rencontr\'{e} des difficult\'{e}s \`{a} r\'{e}gler le correcteur de la
boucle externe car seule une m\'{e}thode empirique est applicable. Nous
avons respect\'{e} la r\`{e}gle classique stipulant que la dynamique du
correcteur de la boucle externe est moins rapide que la boucle
interne.

\subsection{Commande \'{e}chantillonn\'{e}e bloqu\'{e}e \`{a} 2 minutes}
\begin{figure*}
\subfigure[D\'ebit amont]{\rotatebox{-0}{\includegraphics*[width=
0.725\columnwidth]{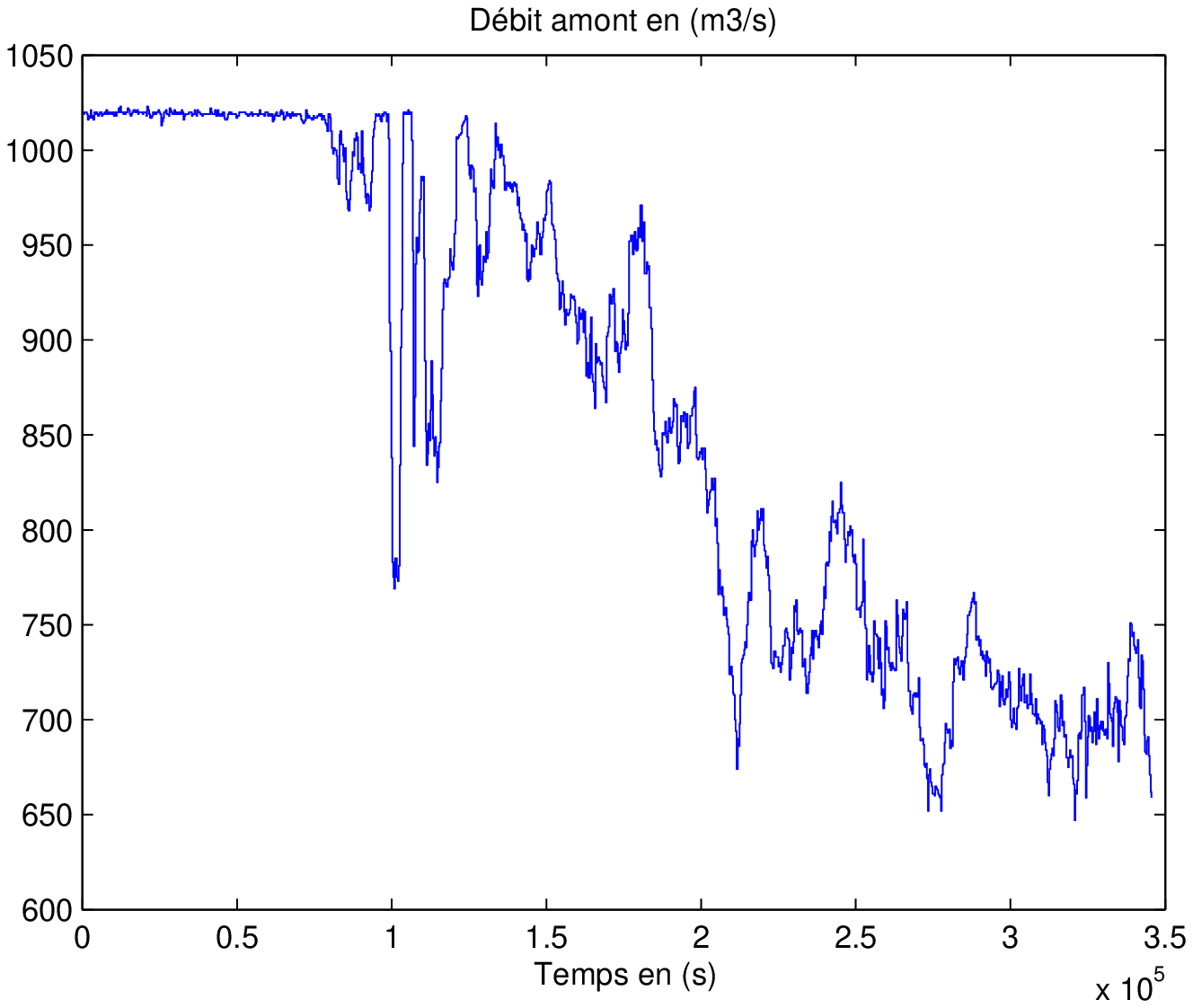}}} \subfigure[Perturbations : biais et
sass\'{e}es]{\rotatebox{-0}{\includegraphics*[width=
0.725\columnwidth]{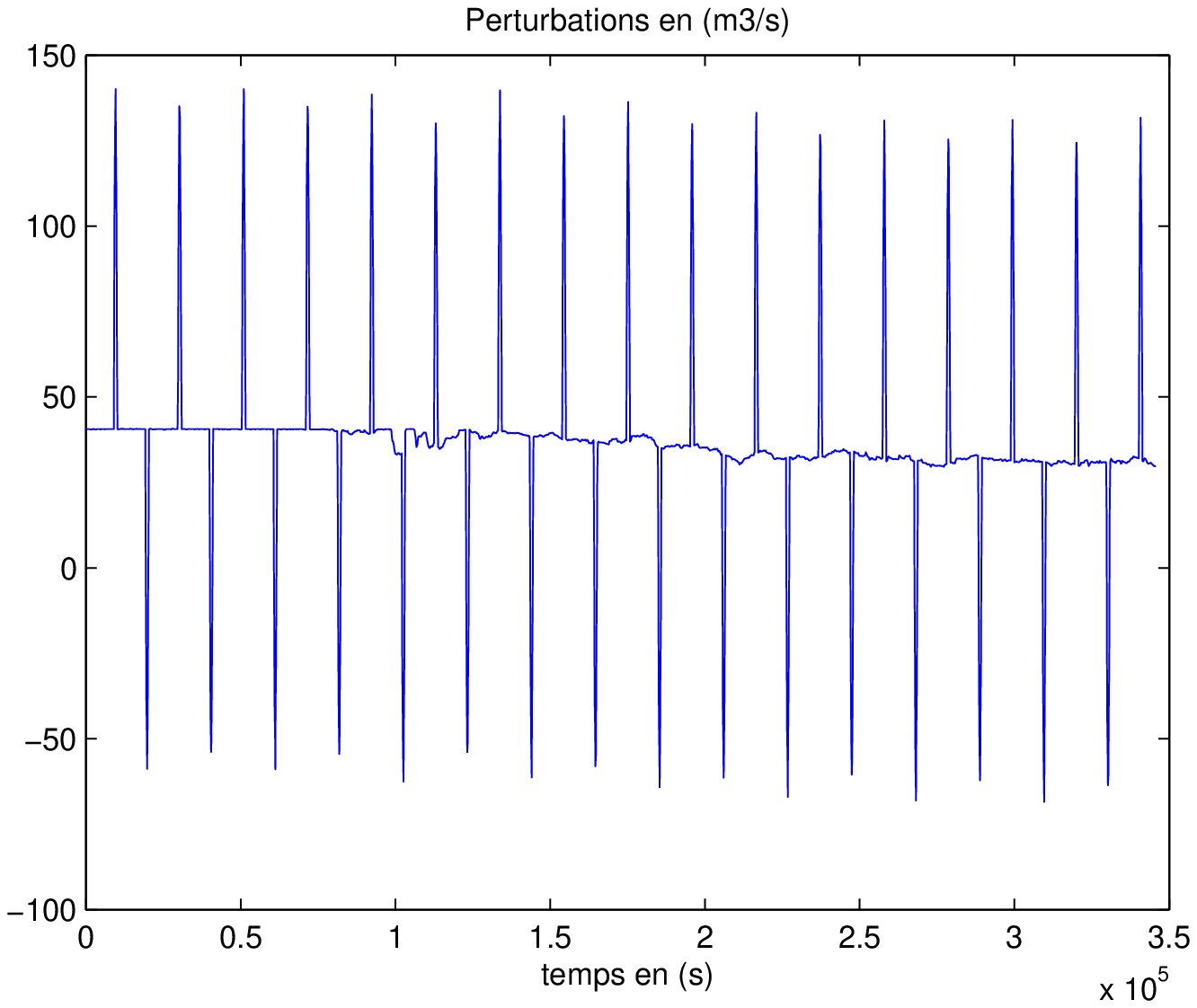}}}
\subfigure[Commande]{\rotatebox{-0}{\includegraphics*[width=
0.725\columnwidth]{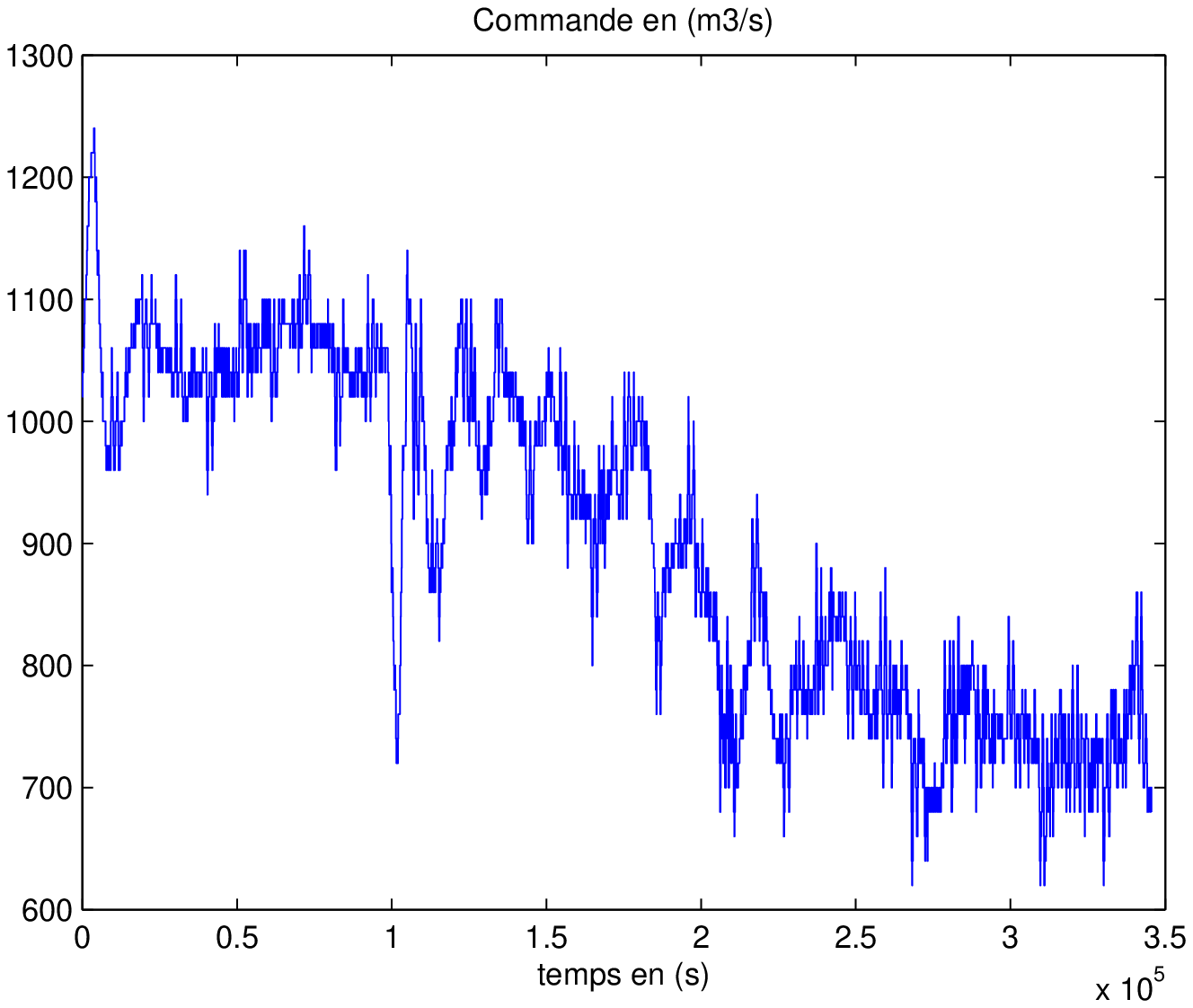}}}\\
\center\subfigure[Poursuite de
$z$]{\rotatebox{-0}{\includegraphics*[width=0.725\columnwidth]{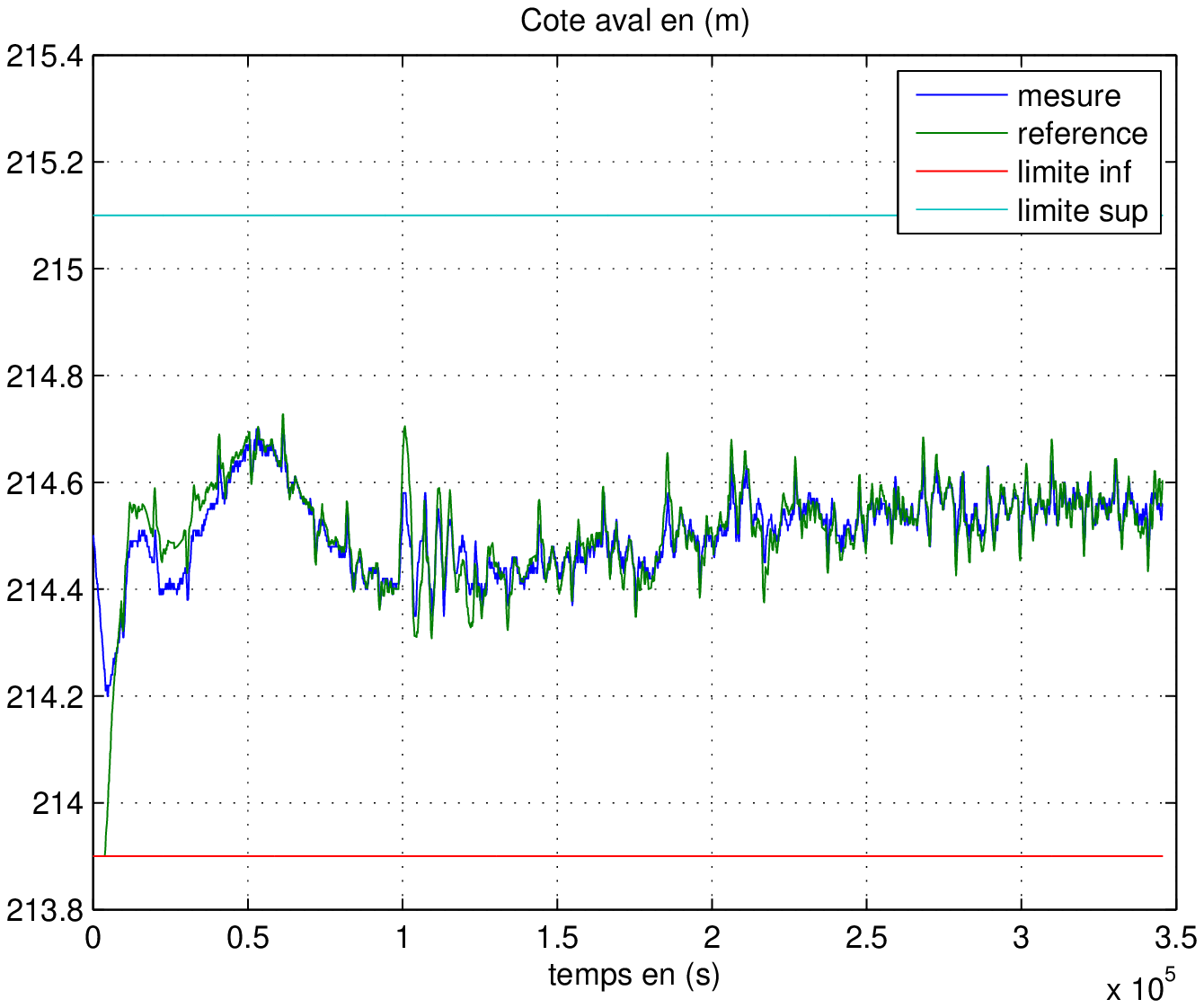}}}\subfigure[Poursuite
de
$z_r$]{\rotatebox{-0}{\includegraphics*[width=0.725\columnwidth]{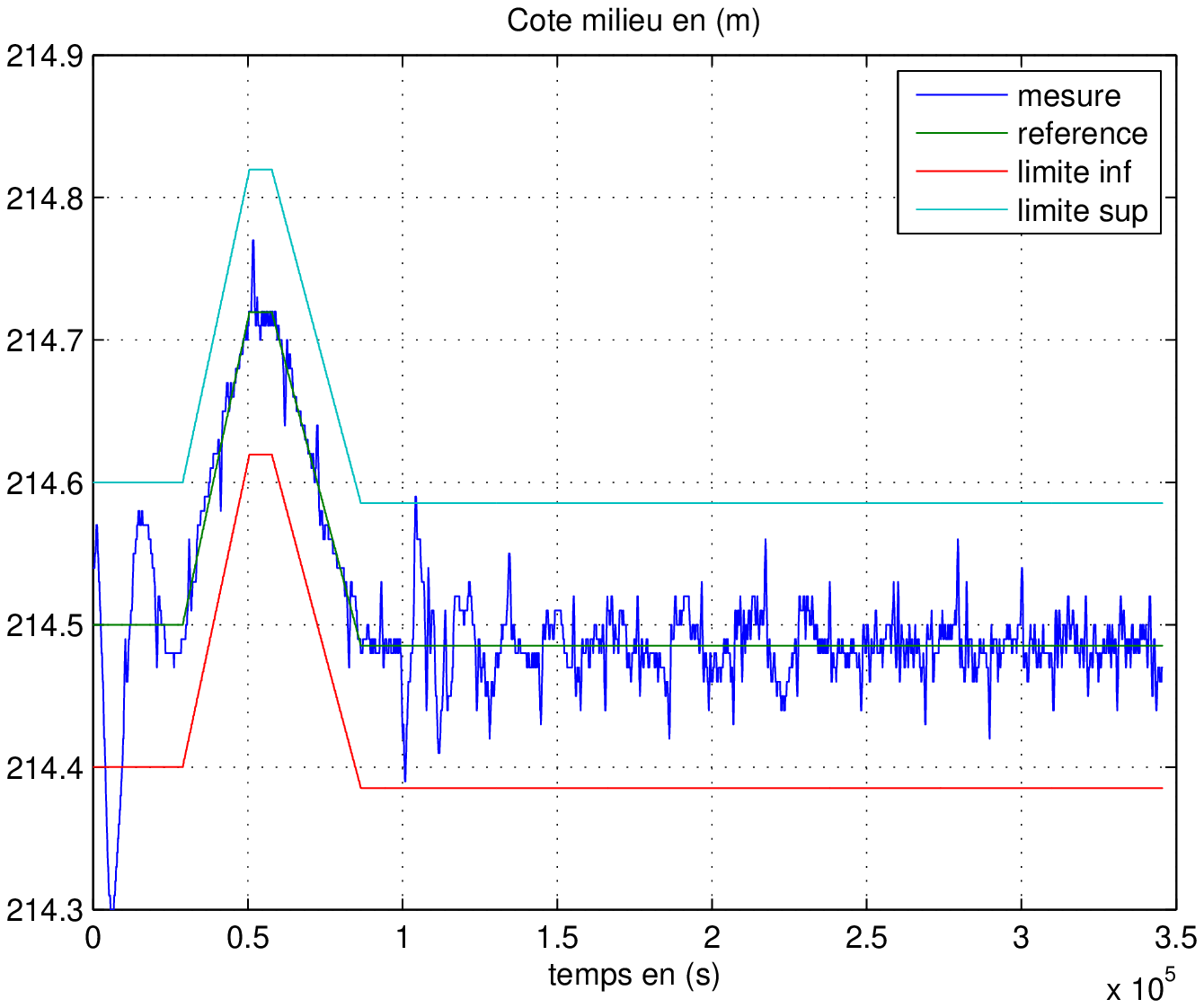}}}
 \caption{Sc\'{e}narion 1 avec poursuite et rejet de fortes perturbations:
 p\'eriode de 2min\label{s13}}
\end{figure*}

\begin{figure*}
\subfigure[D\'ebit amont]{\rotatebox{-0}{\includegraphics*[width=
0.725\columnwidth]{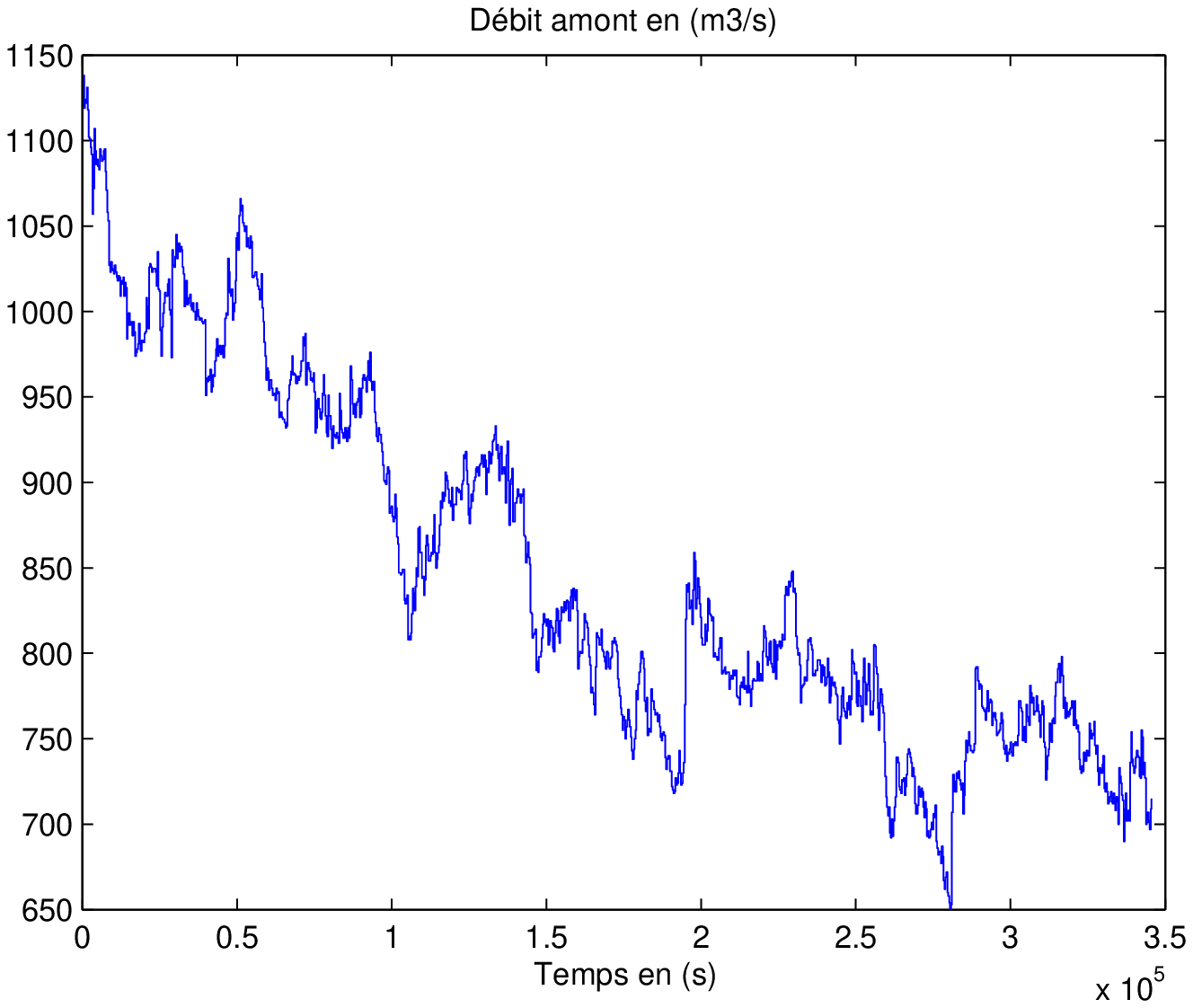}}} \subfigure[Perturbations : biais et
sass\'{e}es]{\rotatebox{-0}{\includegraphics*[width=
0.725\columnwidth]{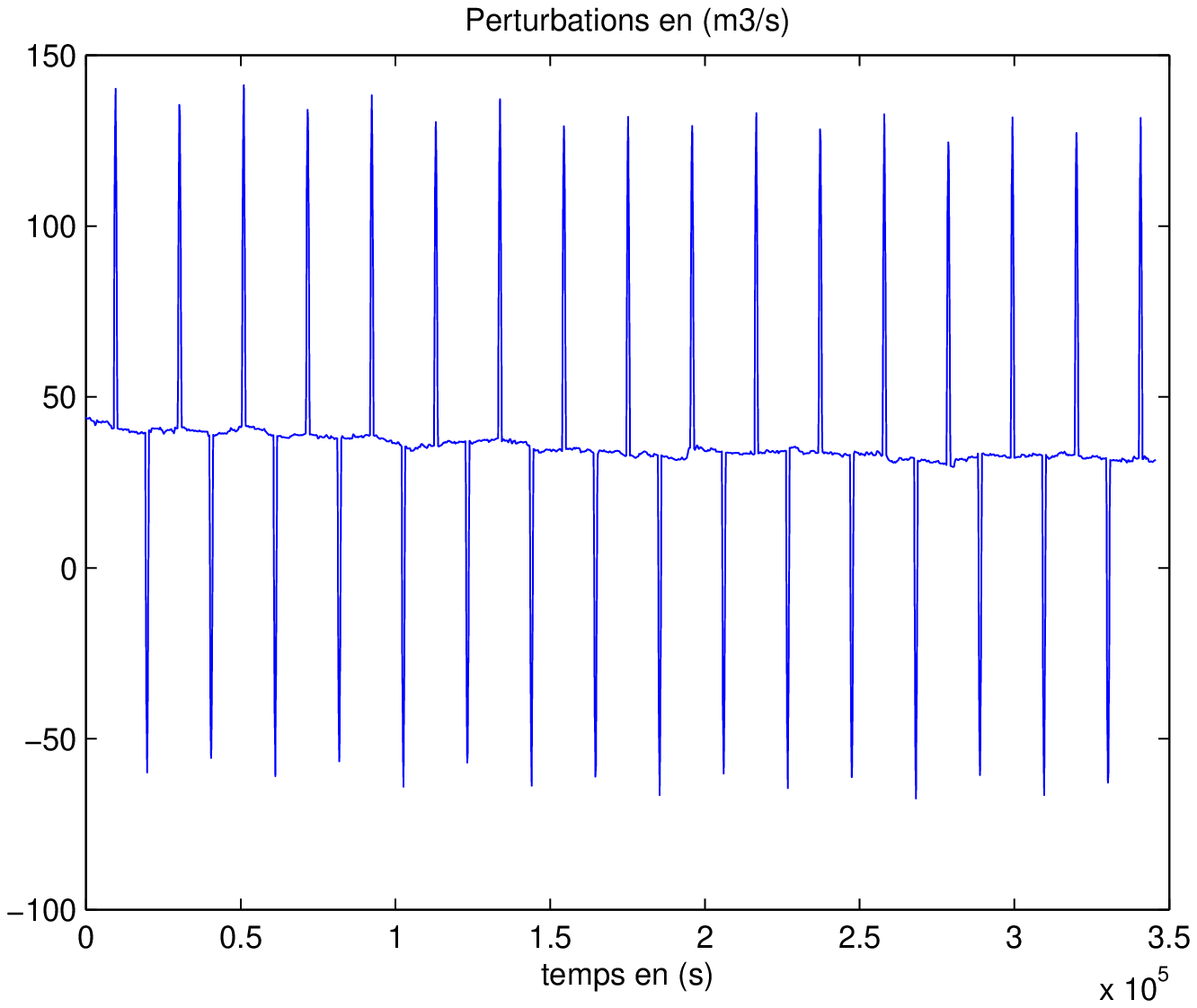}}}
\subfigure[Commande]{\rotatebox{-0}{\includegraphics*[width=
0.725\columnwidth]{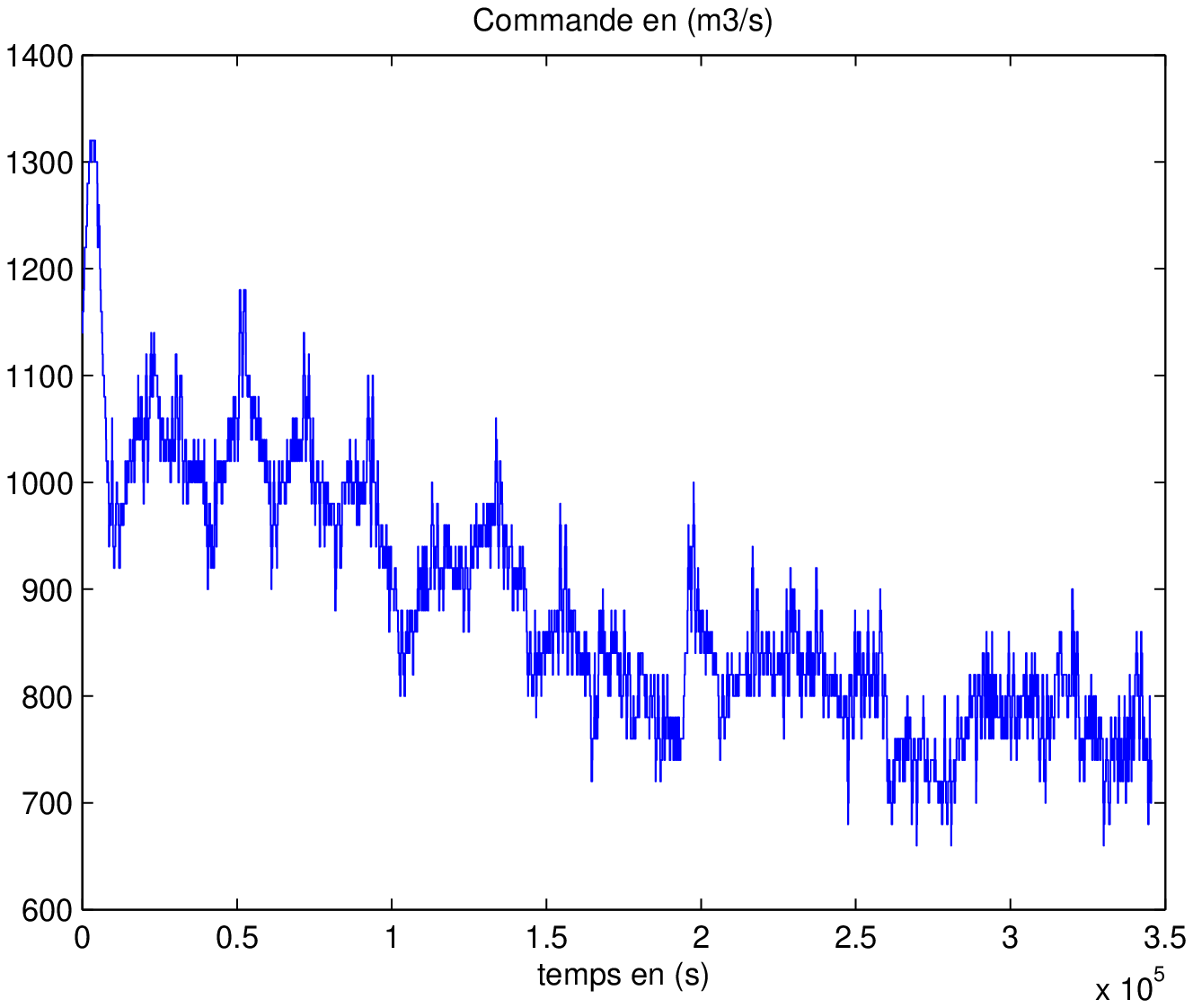}}} \center\subfigure[Poursuite de
$z$]{\rotatebox{-0}{\includegraphics*[width=
0.725\columnwidth]{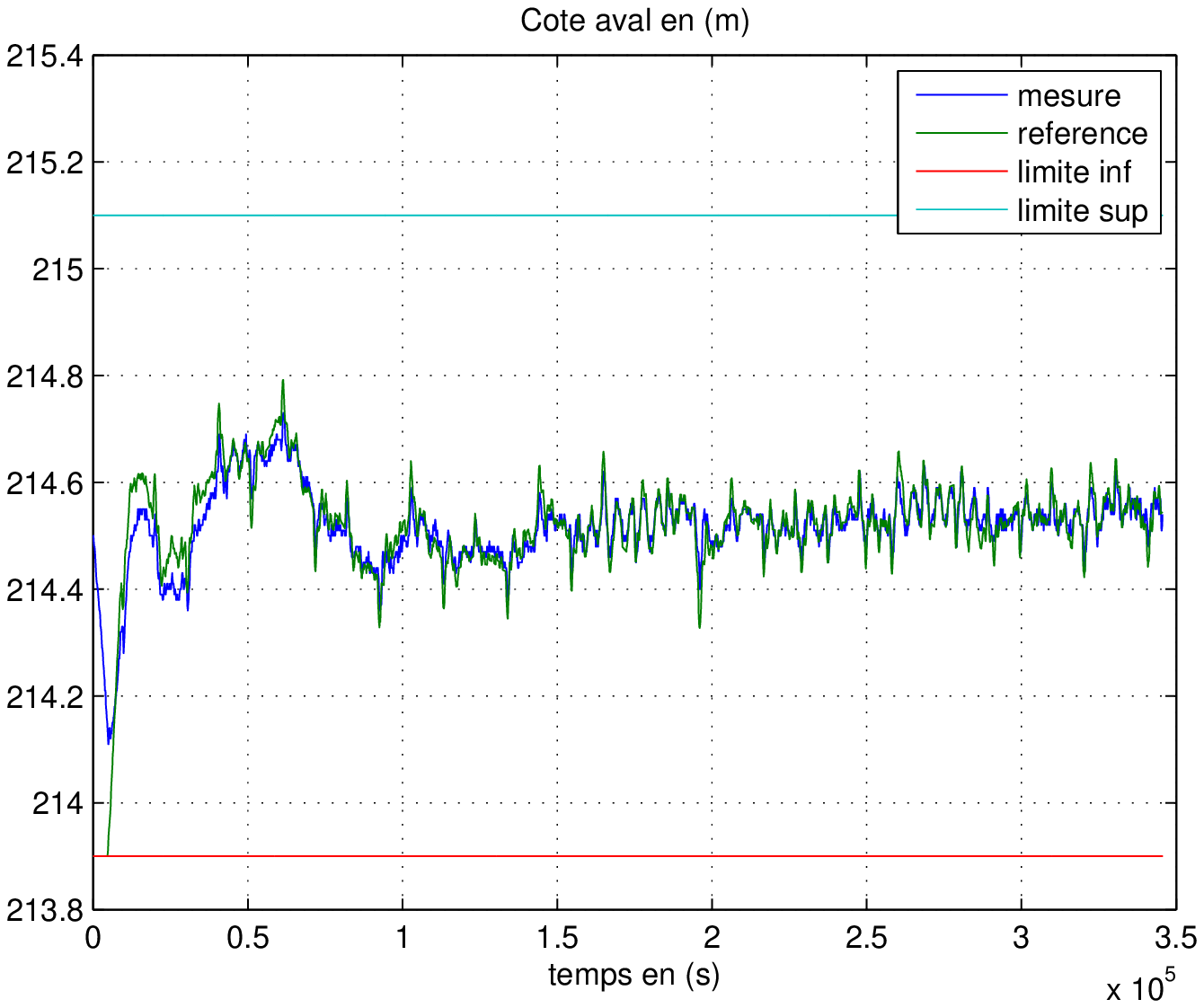}}} \subfigure[Poursuite de
$z_r$]{\rotatebox{-0}{\includegraphics*[width=
0.725\columnwidth]{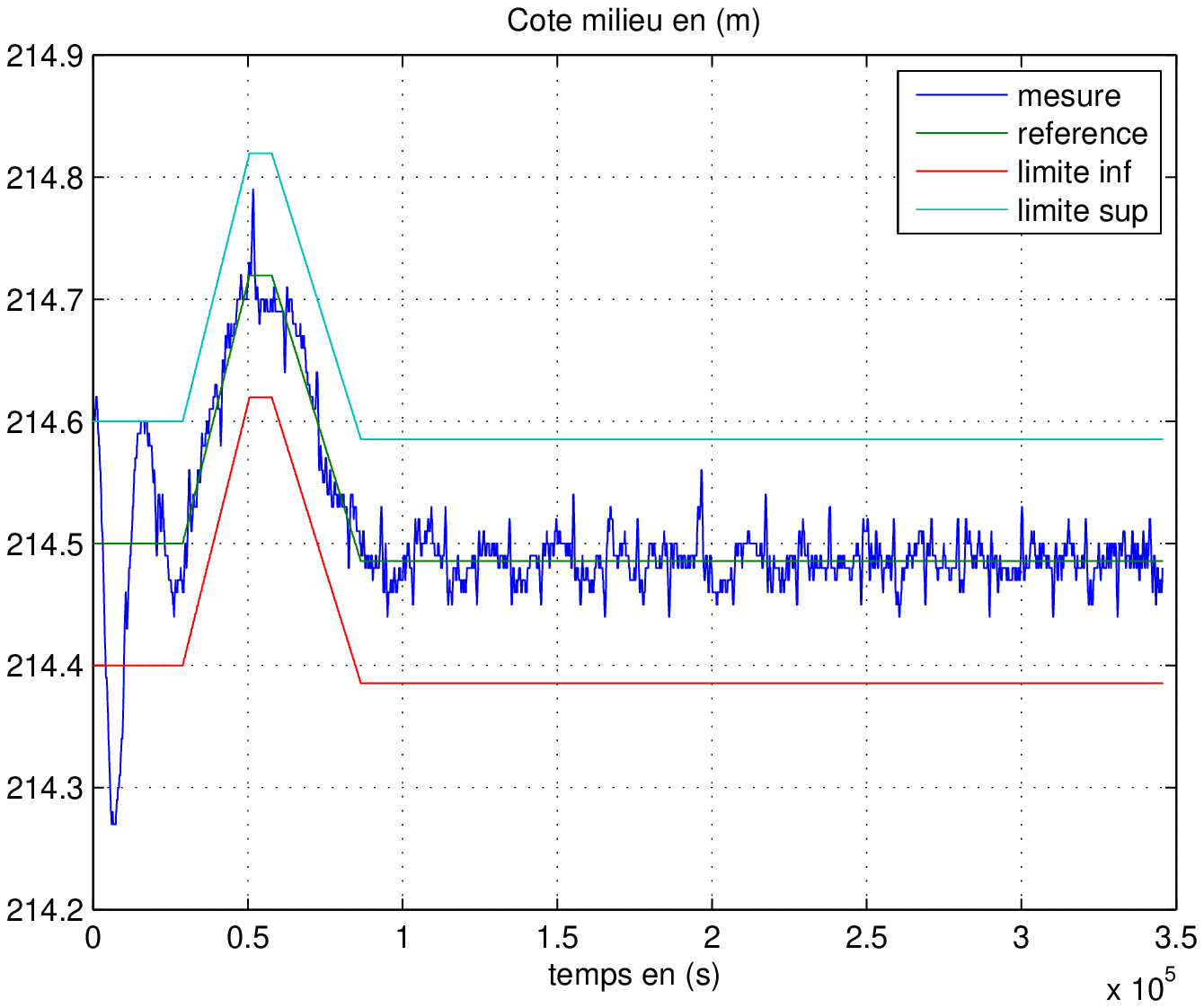}}}
 \caption{Sc\'{e}narion 2 avec poursuite et rejet de perturbations assez
lentes: p\'eriode de 2min\label{s23}}
\end{figure*}

Les Figures \ref{s13} et \ref{s23} pr\'{e}sentent les r\'{e}sultats obtenus
dans le cas o\`{u} la commande est bloqu\'{e}e toutes les 2 minutes, en
respectant les diff\'{e}rents sc\'{e}narios. Aucune contrainte n'est viol\'{e}e:
le contrat est rempli!

\begin{remarque}
Nous avons repris ces m\^{e}mes simulations en augmentant la fr\'{e}quence
d'\'{e}chantillonnage, en passant \`{a} 1 minute. On am\'{e}liore
consid\'{e}rablement \eqref{ineq} en obtenant: $ z_r^\star - 5
~\mbox{\rm cm} < z_r < z_r^\star + 5 ~\mbox{\rm cm} $.
\end{remarque}

\subsection{Anti-emballement}
\begin{figure*}
\subfigure[D\'ebit amont]{\rotatebox{-0}{\includegraphics*[width=
0.725\columnwidth]{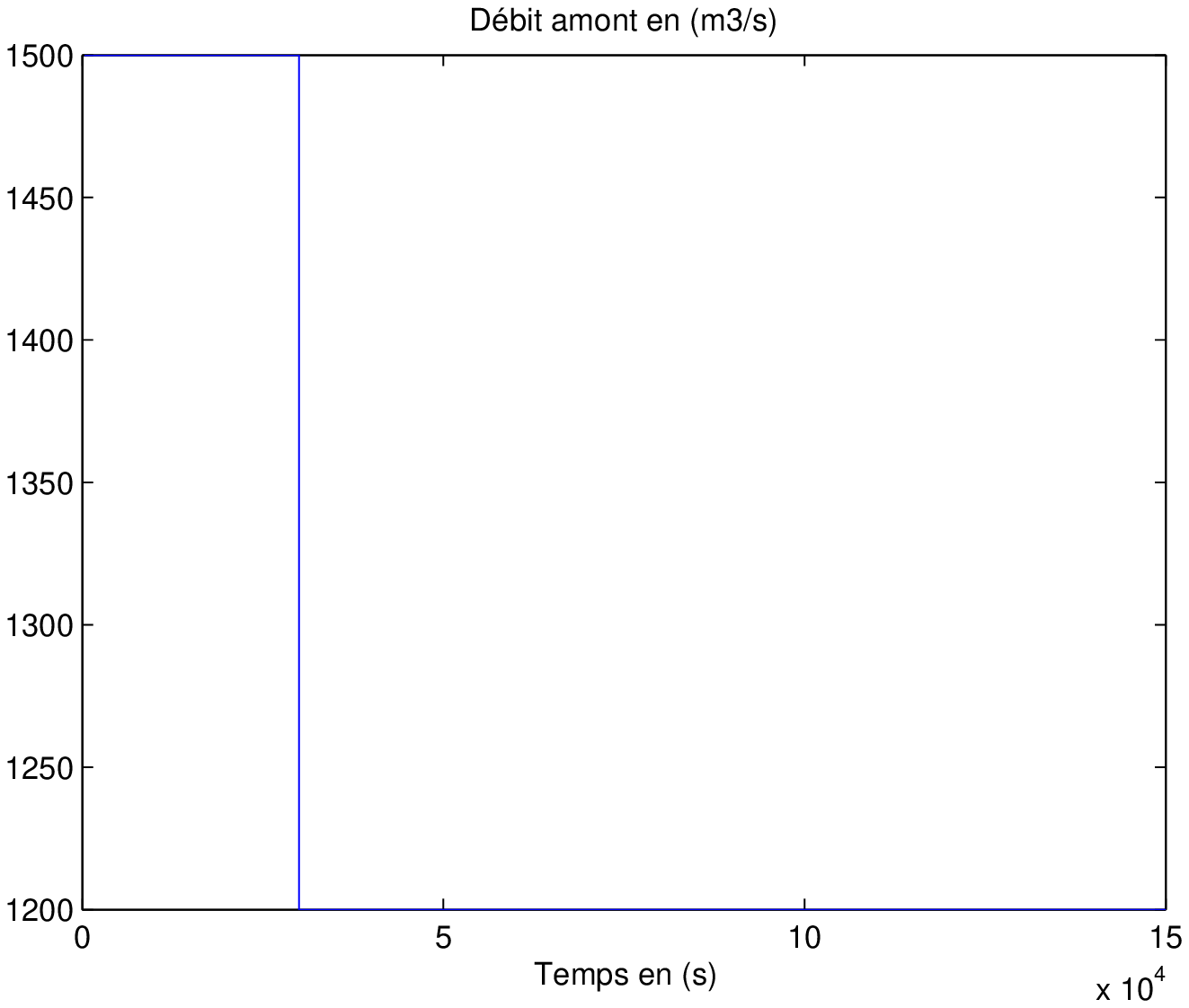}}} \subfigure[Perturbations : biais et
sass\'{e}es]{\rotatebox{-0}{\includegraphics*[width=
0.725\columnwidth]{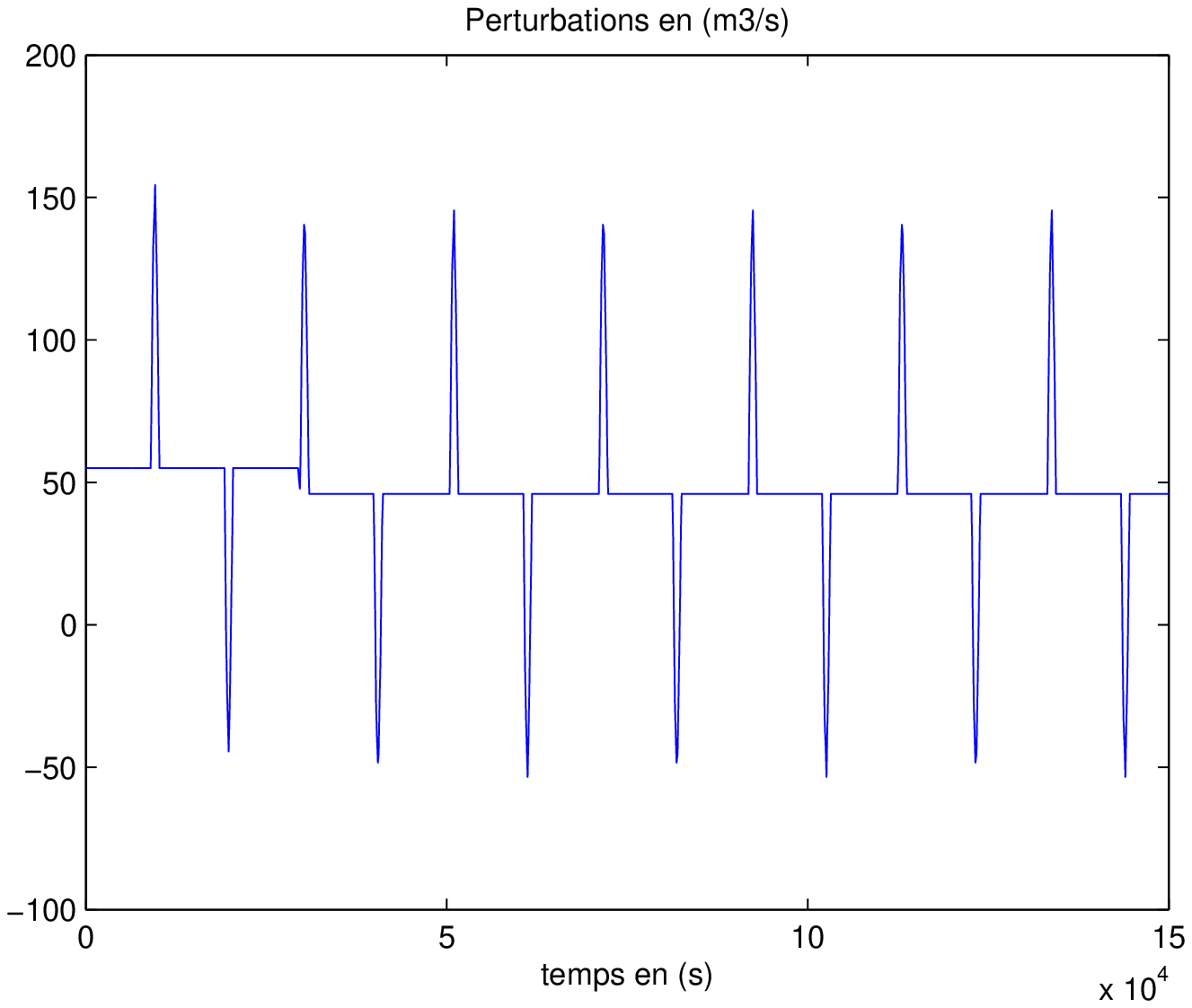}}}
\subfigure[Commande]{\rotatebox{-0}{\includegraphics*[width=
0.725\columnwidth]{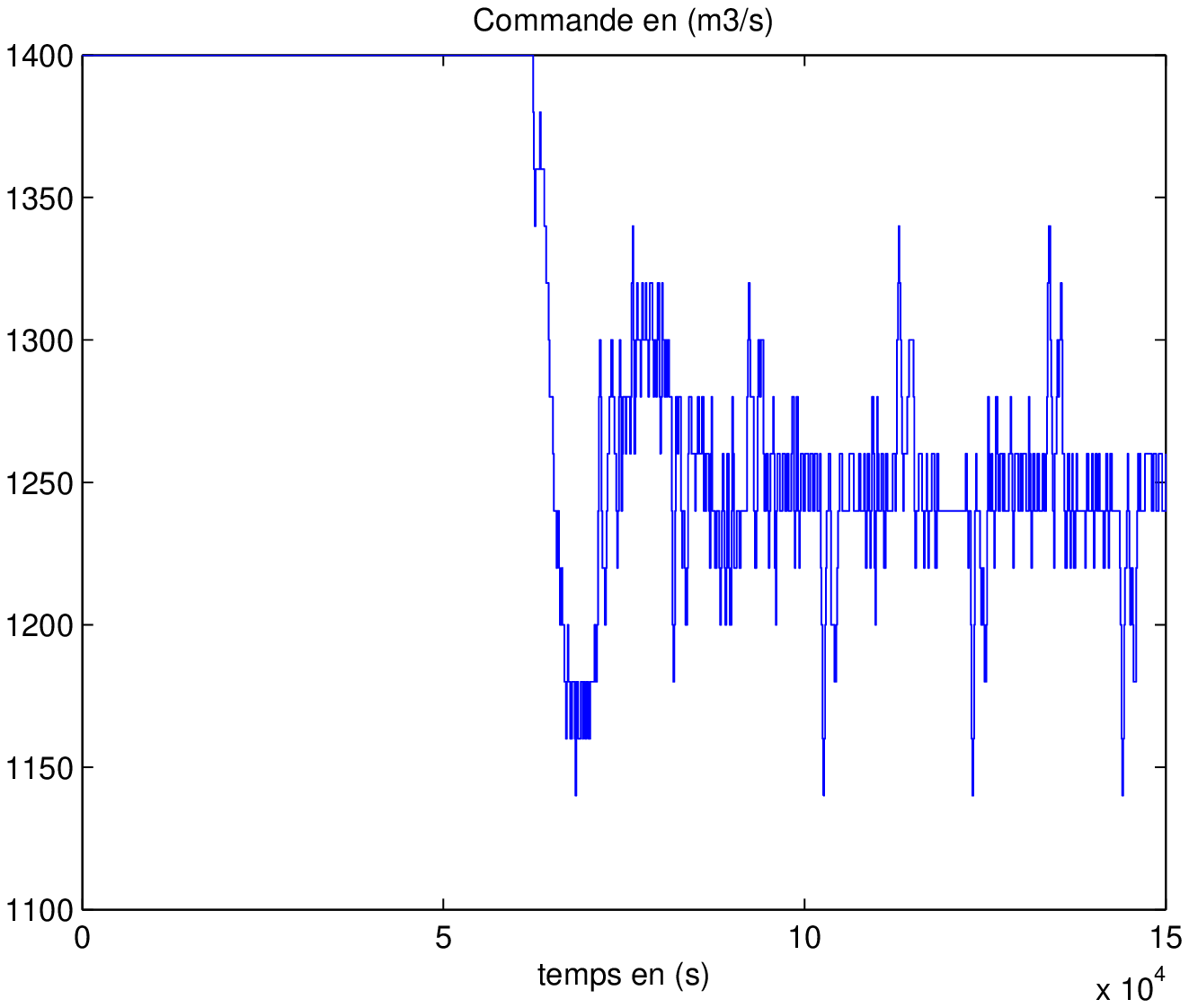}}} \center\subfigure[Poursuite de
$z$]{\rotatebox{-0}{\includegraphics*[width=
0.725\columnwidth]{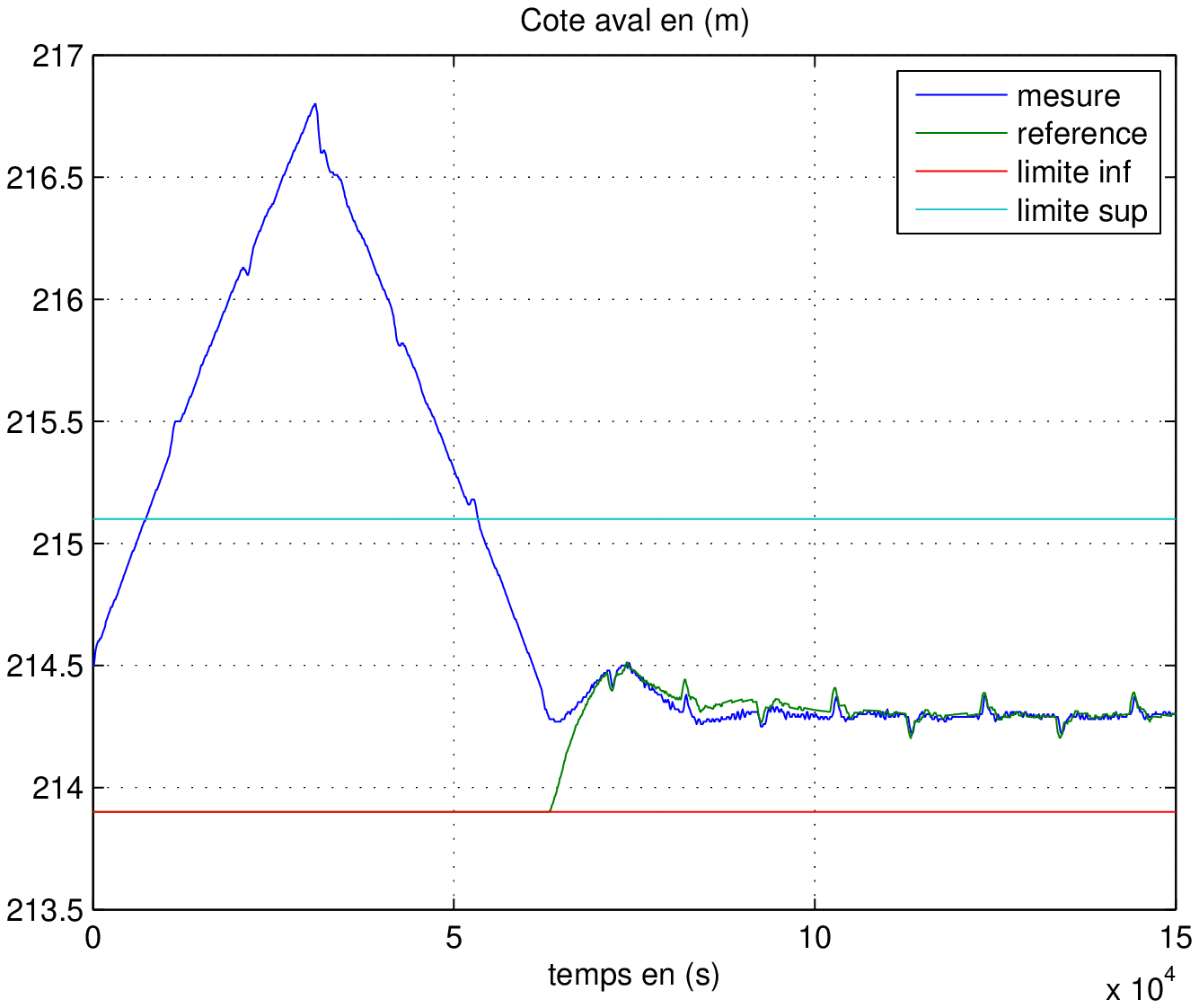}}}\subfigure[Poursuite de
$z_r$]{\rotatebox{-0}{\includegraphics*[width=
0.725\columnwidth]{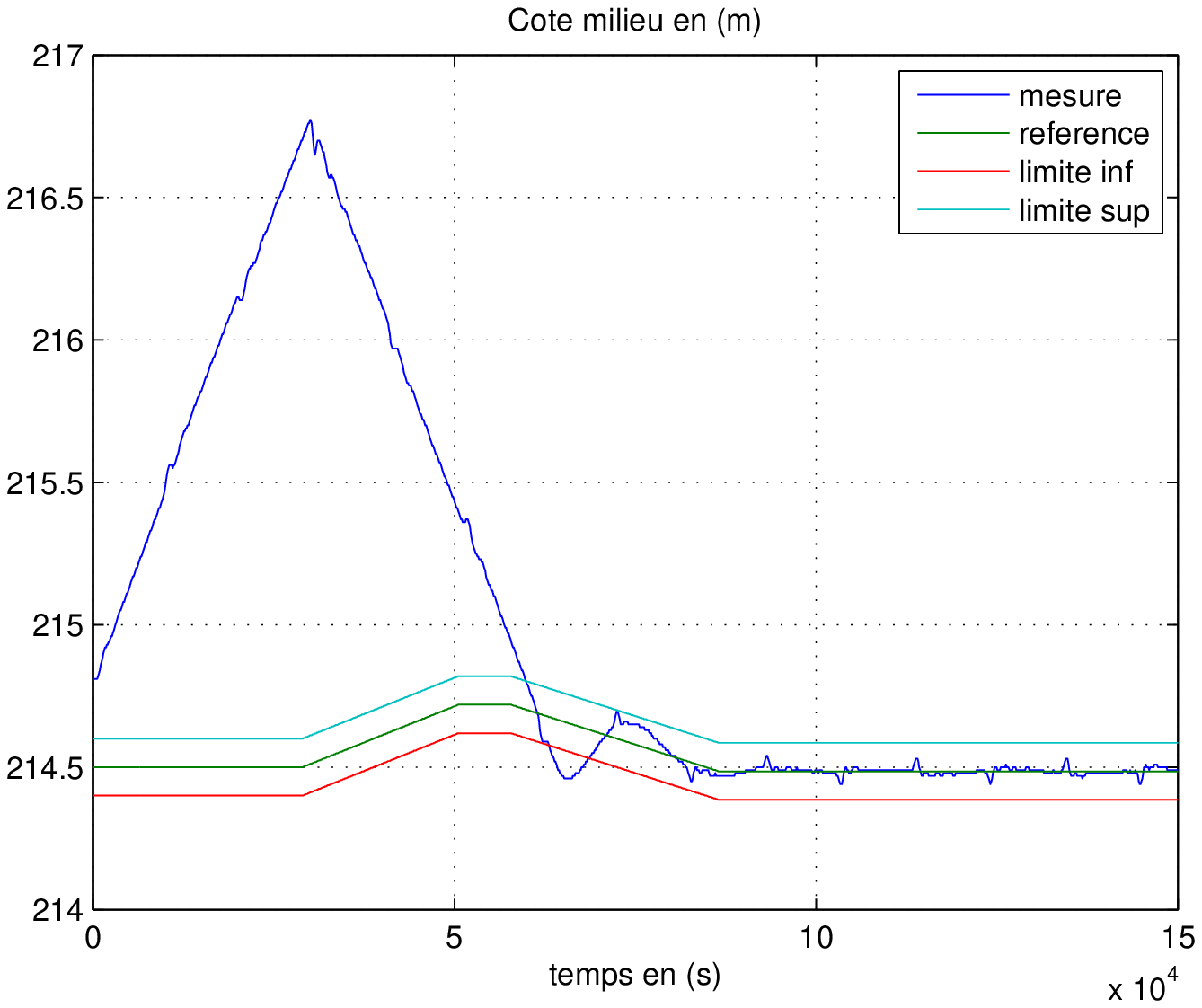}}}
 \caption{Sc\'{e}narion 1 avec poursuite et rejet de perturbations avec actionneur
satur\'{e}: p\'eriode de 2min\label{s33}}
\end{figure*}

Un anti-emballement, ou {\it anti-windup}, est mis en place sur
chaque boucle de r\'{e}gulation (boucles interne et externe). C'est
pourquoi, sur la Figure \ref{s33}, la commande et la replanification
de $z^\ast$ d\'{e}croissent imm\'{e}diatement apr\`{e}s la diminution du d\'{e}bit
d'entr\'{e}e $Q_e$.
%
\subsection{Consigne anticip\'{e}e}
Il est, bien s\^{u}r, possible d'anticiper la consigne de $z_r$: ceci
revient \`{a} pr\'{e}voir les changements de consigne et \`{a} les appliquer,
afin de consid\'{e}rer le temps s\'{e}parant l'\'{e}volution de $z$ de son effet
sur $z_r$. Afin de ne ne pas masquer les r\'{e}sultats de la r\'{e}gulation
propos\'{e}e, cette id\'{e}e n'est pas mise en {\oe}uvre ici. De fait, on voit
appara\^{\i}tre, dans la commande \`{a} $2$ minutes, un retard d'environ $45$
minutes que nous pourrions largement att\'{e}nuer.

\section{Conclusion}\label{conclusion}
\subsection{Int\'{e}r\^{e}t pratique de la commande sans mod\`{e}le}
La strat\'{e}gie propos\'{e}e poss\`{e}de de nombreux avantages:
\begin{itemize}
\item Le retard, cause d'instabilit\'{e} dans bien des lois de commande,
est \'{e}cart\'{e} gr\^{a}ce \`{a} des redondances analytiques approximatives,
d\'{e}duites de la physique du bief. Les erreurs, qui en r\'{e}sultent, sont
corrig\'{e}es par une boucle externe de r\'{e}gulation. Une nouvelle
trajectoire de r\'{e}f\'{e}rence est calcul\'{e}e \`{a} chaque instant pour le
niveau aval.
\item On tire toujours profit de la grande
r\'{e}activit\'{e} de la commande sans mod\`{e}le dans la boucle interne, qui
assure ainsi la poursuite de la trajectoire sans conna\^{\i}tre le mod\`{e}le
du bief.
\end{itemize}
Par ses propri\'{e}t\'{e}s de robustesse, d'adaptabilit\'{e} et de simplicit\'{e},
la commande sans mod\`{e}le apporte des performances remarquables avec
un temps de mise au point tr\`{e}s court, compar\'{e} aux syst\`{e}mes avanc\'{e}s
de commande, que ce soit en \'{e}tude de simulation ou sur site. Elle
semble donc particuli\`{e}rement adapt\'{e}e au milieu industriel.

L'algorithme de commande propos\'{e} est innovant dans le domaine de la
r\'{e}gulation de niveau des canaux d\'{e}couverts. Les performances
obtenues en poursuite, en rejet de perturbations, et en robustesse
sont remarquables, compte tenu de la s\'{e}v\'{e}rit\'{e} des sc\'{e}narios de
simulation appliqu\'{e}s. Le fait de maintenir le niveau au point milieu
du bief  - soit \`{a} 7,5 kilom\`{e}tres de l'actionneur - dans une bande de
$\pm 10$ centim\`{e}tres, avec des perturbations inconnues, rel\`{e}ve en
effet d'une vraie prouesse. La commande sans mod\`{e}le appliqu\'{e}e \`{a} la
r\'{e}gulation de niveau constitue une r\'{e}gulation industrielle non
seulement en raison des performances atteintes mais aussi, et
surtout, pour ses qualit\'{e}s intrins\`{e}ques:
\begin{itemize}
\item facilit\'{e} de mise en oeuvre (structure de type PID
{\it feedforward}),

\item faible sollicitation de la charge du calculateur,

\item temps de mise au point r\'{e}duit, gr\^{a}ce, notamment, au faible nombre de
param\`{e}tres \`{a} r\'{e}gler\footnote{5 si l'on fait la somme des param\`{e}tres
dans \eqref{F} et \eqref{universal}. Si l'on part de \eqref{sm}, on
n'a plus, comme ici, que 3 param\`{e}tres \`{a} choisir. Le calibrage
devient trivial car les gains de l'i-PI servent, d'apr\`{e}s la remarque
\ref{trivial}, \`{a} r\'{e}guler un int\'{e}grateur pur du premier ordre.},

\item auto-adaptation par rapport aux variations du proc\'{e}d\'{e},

\item maintenance ais\'{e}e, gr\^{a}ce \`{a} un algorithme tr\`{e}s simple.
\end{itemize}
L'int\'{e}r\^{e}t pour le producteur hydraulique EDF est de disposer d'une
commande capable de r\'{e}guler un niveau proche ou distant du vannage
(organe r\'{e}glant), tout en garantissant les performances d\'{e}sir\'{e}es sur
toute la plage de fonctionnement de l'usine hydro\'{e}lectrique, malgr\'{e}
la pr\'{e}sence de perturbations impr\'{e}visibles. Ces atouts permettront
certainement de r\'{e}duire les co\^{u}ts des projets de construction ou de
r\'{e}novation d'installations hydrauliques, en diminuant la dur\'{e}e des
essais de mise en service.

\subsection{Quel r\^{o}le pour la dimension infinie?}
Des techniques \'{e}l\'{e}mentaires suffisent, ici, pour traiter d'un sujet
qui rel\`{e}ve de la dimension infinie (syst\`{e}mes \`{a} retards et r\'{e}gis par
des \'{e}quations aux d\'{e}riv\'{e}es partielles). Reste \`{a} savoir si cette
constatation peut s'\'{e}tendre \`{a} d'autres applications dont la
mod\'{e}lisation naturelle requiert {\it a priori} la dimension infinie.
C'est, sans aucun doute, une question m\'{e}thodologique, voire
\'{e}pist\'{e}mologique, fondamentale pour le futur de l'automatique et,
peut-\^{e}tre, de plusieurs autres domaines d'ing\'{e}nierie et de
math\'{e}matiques appliqu\'{e}es. Seul l'examen de multiples cas concrets,
pos\'{e}s par les praticiens, permettra d'y r\'{e}pondre.


\end{document}